\documentclass[pdflatex,sn-mathphys-num]{sn-jnl}

\usepackage{graphicx}%
\usepackage{multirow}%
\usepackage{amsmath,amssymb,amsfonts}%
\usepackage{amsthm}%
\usepackage{mathrsfs}%
\usepackage[title]{appendix}%
\usepackage{xcolor}%
\usepackage{textcomp}%
\usepackage{manyfoot}%
\usepackage{booktabs}%
\usepackage{algorithm}%
\usepackage{algorithmicx}%
\usepackage{algpseudocode}%
\usepackage{listings}%

\theoremstyle{thmstyleone}%
\newtheorem{theorem}{Theorem}
\newtheorem{lemma}[theorem]{Lemma}%
\newtheorem{corollary}[theorem]{Corollary}%

\theoremstyle{thmstyletwo}%

\theoremstyle{thmstylethree}%

\newcommand{\R}{\mathbb{R}}

\newcommand{\Rnn}{\mathbb{R}^{n}_{+}}

\newcommand{\vDC}{v_{\mathrm{DC}}}
\newcommand{\vE}{v_{\mathcal{E}}}

\let\papertextv\v
\DeclareRobustCommand{\v}{%
  \ifmmode
    v%
  \else
    \expandafter\papertextv
  \fi
}
\DeclareMathOperator{\tr}{tr}

\DeclareMathOperator{\diag}{diag}

\usepackage{xspace}

\usepackage[nameinlink,capitalise]{cleveref}

\NewCommandCopy{\papercite}{\cite}
\RenewDocumentCommand{\cite}{o m}{%
  \IfNoValueTF{#1}{\papercite{#2}}{\papercite[][#1]{#2}}%
}
\let\savedproof\proof
\let\savedendproof\endproof
\AtBeginDocument{%
  \let\proof\savedproof
  \let\endproof\savedendproof
}

\begin{document}

\title{The Optimal D.C. Bound Is Dominated by
the Reciprocal of the Lov\'asz Theta Number}

\author[1]{\fnm{Yuqin} \sur{Ao}}
\author[1]{\fnm{Yong} \sur{Xia}}
\affil[1]{\orgdiv{School of Mathematical Sciences},
  \orgname{Beihang University}, \orgaddress{\city{Beijing}, \country{China}}}
\abstract{We resolve a question left open by Anstreicher and Burer
(J. Global Optim. 33, 299--312, 2005). For the standard quadratic
formulation of the stability number, we prove that the optimal
difference-of-convex (D.C.) bound is at most the reciprocal of the
Lov\'asz theta number. The dominance is strict whenever the
Lov\'asz theta number is not exact.}

\maketitle

\section{Introduction}\label{sec:intro}
For a symmetric matrix $Q\in\R^{n\times n}$, the standard
quadratic program is
\begin{equation}\label{eq:intro-qps}
  v_{\mathrm{QPS}}(Q)
  =\min\{x^TQx:~x\in\Delta\},
  \qquad
  \Delta:=\{x\in\Rnn:e^Tx=1\},
\end{equation}
where $\Rnn$ denotes the nonnegative orthant and $e$ is the
all-ones vector.

Any quadratic program over a nonempty bounded polyhedron
can be equivalently reformulated in this form by expressing
each feasible point as a convex combination of the vertices.
An example is quadratic optimization over the
$\ell_1$-ball~\cite{bomze2007,PT2006,xia2015}.

In this paper, we consider the graph case $Q=I+A$.
Let $G=(V,\mathcal{E})$ be a finite simple graph with
$V=\{1,\ldots,n\}$, $n\geq1$, adjacency matrix $A$, and
stability number $\alpha(G)$, the maximum cardinality of
a set of pairwise nonadjacent vertices.
The Motzkin--Straus theorem~\cite{motzkin1965} gives
\begin{equation}\label{eq:intro-motzkin-straus}
  v_{\mathrm{QPS}}(I+A)=\frac{1}{\alpha(G)}.
\end{equation}
Thus, any positive lower bound $L$ on
$v_{\mathrm{QPS}}(I+A)$ yields the upper bound
$\alpha(G)\leq1/L$.

The clique cover number $\bar\chi(G)$ is the minimum number of
cliques needed to cover all vertices of $G$.
Its fractional relaxation $\bar\chi_f(G)$ minimizes the total
weight assigned to cliques, subject to nonnegative weights and
a total weight of at least one on the cliques containing each vertex.
These bounds satisfy
\[
  \alpha(G)\leq\bar\chi_f(G)\leq\bar\chi(G),
\]
since each clique contains at most one vertex of a stable set
and any clique cover yields a feasible fractional cover by
assigning unit weight to each selected clique.
Computing $\bar\chi(G)$ is NP-hard because it is the chromatic
number of the complement graph~\cite{karp1972}. Computing $\bar\chi_f(G)$ is also NP-hard because it is
the fractional chromatic number of the complement
graph~\cite{grotschel1981}.

\textbf{Some notation.}
 $X\succeq0$ denotes positive semidefiniteness and
$X\geq0$ denotes entrywise nonnegativity. The matrix inner
product is written $A\bullet B=\operatorname{tr}(AB)$, where
$\operatorname{tr}(\cdot)$ denotes the trace. Set $E=ee^T$. The operators $\diag$
and $\operatorname{Diag}$ extract a matrix diagonal and form a diagonal
matrix from a vector, respectively. For $S\subseteq V$, $G[S]$ is the
induced subgraph, $M_{SS}$ is the principal submatrix indexed by $S$, and
$v_S$ is the corresponding subvector of $v$. We write $P_3$ for the path
on three vertices and $u^\perp$ for the orthogonal complement of $u$.

The Lov\'asz theta number~\cite{lovasz1979} is an upper
bound on $\alpha(G)$ that can be computed to any prescribed
accuracy in polynomial time~\cite[Section~6]{grotschel1981}.
It is at least as strong as the fractional clique cover
bound~\cite[Theorem~10]{lovasz1979}:
\[
  \alpha(G)\leq\vartheta(G)\leq\bar\chi_f(G).
\]
It admits the following semidefinite programming formulation:
\begin{equation}\label{eq:intro-theta}
  \begin{aligned}
    \vartheta(G)=\max_X\quad &E\bullet X\\
    \text{s.t.}\quad &\operatorname{tr}X=1,\\
                    &X_{ij}=0,\qquad \{i,j\}\in \mathcal{E},\\
                    &X\succeq 0.
  \end{aligned}
\end{equation}

Adding $X\geq0$ to \eqref{eq:intro-theta} yields Schrijver's
bound $\vartheta'(G)$~\cite{schrijver1979}, with
$\alpha(G)\leq\vartheta'(G)\leq\vartheta(G)$.

 Let $\mathcal{CP}_n=\operatorname{cone}\{xx^T:x\in\Rnn\}$. The completely
positive reformulation of \eqref{eq:intro-qps} is~\cite{bomze2000}
\[
  v_{\mathrm{QPS}}(Q)
  = \min_X\{Q\bullet X:E\bullet X=1,\ X\in\mathcal{CP}_n\}.
\]
Since every completely positive matrix is positive semidefinite and
entrywise nonnegative, the doubly nonnegative (DNN) relaxation gives
\[
  v_{\mathrm{DNN}}(Q)
  :=\min_X\{Q\bullet X:E\bullet X=1,\ X\succeq0,\ X\geq0\}\leq v_{\mathrm{QPS}}(Q).
\]

Xia et al.~\cite{xia2013} estimate the gap
$v_{\mathrm{QPS}}(Q)-v_{\mathrm{DNN}}(Q)$ from below to improve the
bound, and G\"okmen and Y{\i}ld{\i}r{\i}m~\cite{gokmenYildirim2022}
characterize when the gap vanishes. The DNN bound can also be viewed as
the minimum of a convex underestimator on the simplex~\cite{yildirim2022}.

For the stability-number formulation, de Klerk and
Pasechnik~\cite[Remark~5.1]{deKlerkPasechnik2002} give
\[
  v_{\mathrm{DNN}}(G):=v_{\mathrm{DNN}}(I+A)=1/\vartheta'(G).
\]

Convex quadratic underestimators also play a role in global
solution methods~\cite{judiceSessaFukushima2026}.

Bomze~\cite{bomze2002} constructed lower bounds for standard
quadratic programs by decomposing $Q=S-T$, where $S,T\succeq0$,
and minimizing the convex and concave terms separately.
Since the minimum of $-x^TTx$ over $\Delta$ is attained at
a vertex, this yields
\[
  v(S,T)
  :=\min_{x\in\Delta}x^TSx-\max_iT_{ii}
  \leq v_{\mathrm{QPS}}(Q).
\]

Anstreicher and Burer~\cite{anstreicher2005} showed that the
optimal bound over all such decompositions,
\[
  v_{\mathrm{DC}}(Q)
  :=\sup\{v(S,T):Q=S-T,\ S,T\succeq0\},
\]
admits a semidefinite programming formulation.
For $Q=I+A$, write $v_{\mathrm{DC}}(G):=v_{\mathrm{DC}}(I+A)$.
Their single-matrix formulation gives
\begin{equation}\label{eq:intro-dc}
  \begin{aligned}
    v_{\mathrm{DC}}(G)=\min_X\quad &(I+A)\bullet X\\
    \text{s.t.}\quad &E\bullet X=1,\\
                    &X\succeq0,\quad Xe\geq0,\\
                    &\operatorname{Diag}(Xe)-X\succeq0.
  \end{aligned}
\end{equation}

They also established that
$v_{\mathrm{DC}}(G)\leq v_{\mathrm{DNN}}(G)=1/\vartheta'(G)$.
To compare the optimal D.C. bound with the reciprocal
of the Lov\'asz theta number, Anstreicher and
Burer~\cite[Section~5]{anstreicher2005} introduced
the edge relaxation
\begin{equation}\label{eq:intro-edge}
  \begin{aligned}
    v_{\mathcal E}(G)=\min_X\quad &(I+A)\bullet X\\
    \text{s.t.}\quad &E\bullet X=1,\\
                    &X\succeq0,\\
                    &X_{ij}\geq0,\qquad \{i,j\}\in\mathcal E,
  \end{aligned}
\end{equation}
and proved that
$v_{\mathcal E}(G)=1/\vartheta(G)$.
They reported finding no graph with
$v_{\mathrm{DC}}(G)>v_{\mathcal E}(G)$, but were unable
to prove that $v_{\mathrm{DC}}(G)\leq v_{\mathcal E}(G)$
holds for every graph.

In this paper, we resolve this question by proving that,
for every finite simple graph $G$ with at least one vertex,
\begin{equation}\label{eq:intro-comparison}
  v_{\mathrm{DC}}(G)\leq v_{\mathcal E}(G)
  =\frac{1}{\vartheta(G)}.
\end{equation}
Moreover, we establish strict dominance whenever
the Lov\'asz bound is not exact:
\begin{equation}\label{eq:intro-strict}
  \vartheta(G)>\alpha(G)
  \quad\Longrightarrow\quad
  v_{\mathrm{DC}}(G)<\frac{1}{\vartheta(G)}.
\end{equation}
Consequently, the two bounds coincide if and only if
the D.C. bound is exact.

\section{Dominance and strictness}\label{sec:main}

We first establish the comparison
\eqref{eq:intro-comparison} by constructing a D.C.-feasible
matrix from a theta-optimal solution.
We then use a minimum-support argument and an induced-path
perturbation to prove the strict inequality
\eqref{eq:intro-strict}.

\begin{theorem}\label{thm:main}
For every finite simple graph $G$ with at least one vertex,
\[
  \vDC(G)\leq v_{\mathcal E}(G)=\frac{1}{\vartheta(G)}.
\]
\end{theorem}

\begin{proof}
Write $\vartheta=\vartheta(G)$ and let $Z$ be optimal
in \eqref{eq:intro-theta}. Then
\[
  Z\succeq0,\qquad \tr Z=1,\qquad
  Z_{ij}=0\quad\text{for }\{i,j\}\in\mathcal E,
  \qquad E\bullet Z=\vartheta.
\]
Since $I/n$ is feasible with objective value $1$,
we have $\vartheta\geq1$.
We will show that $X_0:=Z/\vartheta$ is feasible for
the D.C. relaxation \eqref{eq:intro-dc}.

We first derive a consequence of the optimality of $Z$
under diagonal scaling.
For any $a\in\R^n$, set
\[
  D_a=\operatorname{Diag}(a),\qquad Z_a=D_aZD_a.
\]
Then $Z_a\succeq0$ and
$(Z_a)_{ij}=a_i a_jZ_{ij}=0$ on every edge.
If $\tr Z_a>0$, the normalized matrix
$Z_a/\tr Z_a$ is feasible for \eqref{eq:intro-theta}.
Optimality of $Z$ therefore gives
\[
  E\bullet Z_a\leq\vartheta\,\tr Z_a.
\]
If $\tr Z_a=0$, positive semidefiniteness implies $Z_a=0$,
so the same inequality holds.
Since
\[
  E\bullet Z_a=a^TZa,\qquad
  \tr Z_a=a^T\operatorname{Diag}(\diag Z)a,
\]
we obtain
\[
  a^T\bigl(\vartheta\operatorname{Diag}(\diag Z)-Z\bigr)a
  \geq0
  \qquad\text{for every }a\in\R^n.
\]
Consequently,
\begin{equation}\label{eq:theta-slack}
  H:=\vartheta\operatorname{Diag}(\diag Z)-Z\succeq0,
  \qquad
  e^THe=\vartheta\,\tr Z-E\bullet Z=0.
\end{equation}
The latter equality gives
$\|H^{1/2}e\|^2=0$, and hence $He=0$.
Thus
\[
  Ze=\vartheta\diag Z.
\]

We now verify the feasibility of $X_0=Z/\vartheta$.
Positive semidefiniteness and normalization follow from
\[
  X_0\succeq0,\qquad
  E\bullet X_0=\frac{E\bullet Z}{\vartheta}=1.
\]
Moreover, the identity for $Ze$ gives
\[
  X_0e=\diag Z\geq0,
\]
where the inequality follows from $Z\succeq0$, and
\[
  \operatorname{Diag}(X_0e)-X_0
  =\operatorname{Diag}(\diag Z)-\frac{Z}{\vartheta}
  =\frac{H}{\vartheta}\succeq0.
\]
Therefore $X_0$ satisfies all constraints
in \eqref{eq:intro-dc}.

Finally, the edge-zero constraints imply $A\bullet Z=0$.
Evaluating the D.C. objective at $X_0$ yields
\[
  \vDC(G)
  \leq (I+A)\bullet X_0
  =\frac{\tr Z+A\bullet Z}{\vartheta}
  =\frac{1}{\vartheta}=\vE(G).
\]
\end{proof}

We next establish strict dominance.
Choose an optimal solution $Z$ of \eqref{eq:intro-theta}
whose diagonal support $S=\{i:Z_{ii}>0\}$ has minimum
cardinality.
Lemma~\ref{lem:minimum-diagonal-support} shows that
the slack matrix
$H=\vartheta(G)\operatorname{Diag}(\diag Z)-Z$ satisfies
$\ker H_{SS}=\operatorname{span}\{e_S\}$.
This property allows an induced $P_3$ in $G[S]$ to generate
a rank-one perturbation of $X_0=Z/\vartheta(G)$ that preserves
D.C. feasibility and strictly decreases the objective,
as shown in Lemma~\ref{lem:induced-path-gap}.
Finally, Theorem~\ref{thm:strict-dominance} shows that
$\vartheta(G)>\alpha(G)$ guarantees the existence of such
a path, yielding $\vDC(G)<1/\vartheta(G)$.

\begin{lemma}\label{lem:minimum-diagonal-support}
Let $Z$ be optimal in \eqref{eq:intro-theta} with diagonal
support $S=\{i:Z_{ii}>0\}$ of minimum cardinality.
Set
\[
  D=\operatorname{Diag}(\diag Z_{SS}),\qquad
  H_S=\vartheta(G)D-Z_{SS}.
\]
Then $H_S\succeq0$ and
$\ker H_S=\operatorname{span}\{e_S\}$, where $e_S$ is
the all-ones vector indexed by $S$.
\end{lemma}

\begin{proof}
Let $H$ be the matrix in \eqref{eq:theta-slack}
corresponding to $Z$.
Then $H_S=H_{SS}\succeq0$.
Since $Z\succeq0$, we have $Z_{ij}^2\leq Z_{ii}Z_{jj}$.
Thus, if $Z_{ii}=0$, the entire $i$th row and column of $Z$
are zero.
It follows that $H$ vanishes outside $S\times S$.
Together with $He=0$, this gives $H_Se_S=0$.

Suppose that $\ker H_S$ contains a vector
$b\notin\operatorname{span}\{e_S\}$.
Set
\[
  a=b-\left(\min_{i\in S}b_i\right)e_S.
\]
Then $H_Sa=0$, $a\geq0$, and $a$ has both a positive
and a zero entry.
Moreover, $D\succ0$ by the definition of $S$, so
$c:=a^TDa>0$.

Define $\widetilde Z$ by
\[
  \widetilde Z_{SS}
  =\frac{\operatorname{Diag}(a)Z_{SS}\operatorname{Diag}(a)}{c},
  \qquad
  \widetilde Z_{ij}=0
  \quad\text{if }i\notin S\text{ or }j\notin S.
\]
Diagonal scaling and zero extension preserve positive
semidefiniteness and the edge-zero constraints.
Since $H_Sa=0$, we have
$Z_{SS}a=\vartheta(G)Da$, and hence
\[
  \tr\widetilde Z=\frac{a^TDa}{c}=1,\qquad
  E\bullet\widetilde Z
  =\frac{a^TZ_{SS}a}{c}=\vartheta(G).
\]
Therefore $\widetilde Z$ is also optimal
in \eqref{eq:intro-theta}.
However,
\[
  \{i:\widetilde Z_{ii}>0\}
  =\{i\in S:a_i>0\}\subsetneq S,
\]
contradicting the minimum cardinality of $S$.
Thus $\ker H_S=\operatorname{span}\{e_S\}$.
\end{proof}

\begin{lemma}\label{lem:induced-path-gap}
Let $Z$, $S$, and $H_S$ be as in
Lemma~\ref{lem:minimum-diagonal-support}.

Suppose that $i-j-k$ is an induced path in $G[S]$.

Define $a\in\R^S$ by
\[
  (a_i,a_j,a_k)=(1,-2,1),\qquad
  a_\ell=0\quad\text{for }\ell\in S\setminus\{i,j,k\}.
\]
Set $L_S=H_S/\vartheta(G)$ and let $L_S^\dagger$
denote its Moore--Penrose pseudoinverse.
Then
\[
  s:=a^TL_S^\dagger a>0,\qquad
  \vDC(G)\leq\frac{1}{\vartheta(G)}-\frac{2}{s}
  <\frac{1}{\vartheta(G)}.
\]
\end{lemma}

\begin{proof}
Write $\vartheta=\vartheta(G)$ and $Q=I+A$.
Since $i-j-k$ is induced, the only edges between vertices
in the support of $a$ are $\{i,j\}$ and $\{j,k\}$.
Therefore
\[
  e_S^Ta=1-2+1=0,\qquad
  a^TQ_{SS}a
  =a^Ta+2a_i a_j+2a_j a_k
  =6-4-4=-2.
\]

Let $\widehat a\in\R^n$ be the zero extension of $a$,
and consider
\[
  X(t)=X_0+t\widehat a\widehat a^T,\qquad
  X_0=Z/\vartheta,\qquad t\geq0.
\]
By the proof of Theorem~\ref{thm:main}, $X_0$ is feasible
for \eqref{eq:intro-dc} and $Q\bullet X_0=1/\vartheta$.
For every $t\geq0$, we have $X(t)\succeq0$.
Also, $e^T\widehat a=0$ implies
\[
  X(t)e=X_0e\geq0,\qquad
  E\bullet X(t)=E\bullet X_0=1.
\]
It remains to ensure that
$\operatorname{Diag}(X(t)e)-X(t)\succeq0$.

Set $H=\vartheta\operatorname{Diag}(\diag Z)-Z$.
The remaining slack is
\[
  \operatorname{Diag}(X(t)e)-X(t)
  =H/\vartheta-t\widehat a\widehat a^T.
\]
Both terms on the right vanish outside $S\times S$.
Thus this slack is positive semidefinite if and only if
\[
  L_S-taa^T\succeq0.
\]

By Lemma~\ref{lem:minimum-diagonal-support},
$L_S\succeq0$ and
$\ker L_S=\operatorname{span}\{e_S\}$.
Consequently, the range of $L_S$ is $e_S^\perp$.
Since $a\ne0$ and $e_S^Ta=0$, the vector
$b=(L_S^\dagger)^{1/2}a$ satisfies
\[
  L_S^{1/2}b=a,\qquad
  b^Tb=a^TL_S^\dagger a=s>0.
\]
Indeed, $L_S^{1/2}(L_S^\dagger)^{1/2}$ is the orthogonal
projector onto the range of $L_S$.
It follows that
\[
  L_S-\frac{aa^T}{s}
  =L_S^{1/2}
    \left(I-\frac{bb^T}{b^Tb}\right)
    L_S^{1/2}\succeq0,
\]
because the matrix in parentheses is the orthogonal
projector onto $b^\perp$.
Hence $X(1/s)$ is feasible for \eqref{eq:intro-dc}.
Evaluating its objective gives
\[
  \vDC(G)
  \leq Q\bullet X(1/s)
  =Q\bullet X_0+\frac{a^TQ_{SS}a}{s}
  =\frac{1}{\vartheta}-\frac{2}{s}
  <\frac{1}{\vartheta}.\qedhere
\]
\end{proof}

\begin{theorem}\label{thm:strict-dominance}
For every finite simple graph $G$ with at least one vertex,
\[
  \vartheta(G)>\alpha(G)
  \quad\Longrightarrow\quad
  \vDC(G)<\frac{1}{\vartheta(G)}.
\]
\end{theorem}

\begin{proof}
Assume that $\vartheta(G)>\alpha(G)$.
Choose a theta-optimal matrix $Z$ with diagonal support
$S$ of minimum cardinality, as in
Lemma~\ref{lem:minimum-diagonal-support}.
We first show that
$\vartheta(G[S])=\vartheta(G)$.
Since $Z$ vanishes outside $S\times S$, the principal
submatrix $Z_{SS}$ is feasible for the theta program
of $G[S]$ and has objective value $\vartheta(G)$.
Thus $\vartheta(G[S])\geq\vartheta(G)$.
Conversely, extending any theta-feasible matrix for $G[S]$
by zeros gives a theta-feasible matrix for $G$ with
the same objective value.
Hence $\vartheta(G[S])\leq\vartheta(G)$, proving equality.

Suppose that $G[S]$ contains no induced $P_3$.
Then every connected component of $G[S]$ is a clique:
otherwise, a shortest path between two nonadjacent vertices
in a component would have its first three vertices
inducing a $P_3$.
If there are $k$ components, they give a clique cover
of size $k$, while choosing one vertex from each component
gives a stable set of size $k$.
Therefore
\[
  k\leq\alpha(G[S])
  \leq\vartheta(G[S])
  \leq\bar\chi(G[S])
  \leq k.
\]
It follows that
\[
  \vartheta(G)=\vartheta(G[S])
  =\alpha(G[S])\leq\alpha(G),
\]
contradicting the assumption.
Thus $G[S]$ contains an induced $P_3$.
Lemma~\ref{lem:induced-path-gap} now gives
$\vDC(G)<1/\vartheta(G)$.
\end{proof}

\begin{corollary}
For every finite simple graph $G$ with at least one vertex,
\[
  \vDC(G)=\frac{1}{\vartheta(G)}
  \quad\Longleftrightarrow\quad
  \vDC(G)=\frac{1}{\alpha(G)}.
\]
\end{corollary}

\begin{proof}

Suppose that $\vDC(G)=1/\vartheta(G)$.
By Theorem~\ref{thm:strict-dominance},
$\vartheta(G)>\alpha(G)$ is impossible.
Since $\alpha(G)\leq\vartheta(G)$, we obtain
$\vartheta(G)=\alpha(G)$ and hence
$\vDC(G)=1/\alpha(G)$.

Conversely, if $\vDC(G)=1/\alpha(G)$,
Theorem~\ref{thm:main} gives
\[
  \frac{1}{\alpha(G)}
  =\vDC(G)
  \leq\frac{1}{\vartheta(G)}
  \leq\frac{1}{\alpha(G)}.
\]
Equality holds throughout, so
$\vDC(G)=1/\vartheta(G)$.
\end{proof}

\section{Conclusions}\label{sec:con}

We resolve the comparison left open by Anstreicher and Burer:
the optimal D.C. bound is dominated by the reciprocal of the
Lov\'asz theta number, strictly so whenever the Lov\'asz bound
is not exact.
The two bounds coincide if and only if the D.C. bound is exact.
A direction for future research is to identify structural
conditions under which optimized D.C. bounds coincide with,
or are strictly weaker than, semidefinite bounds for more
general nonconvex quadratic programs.

\backmatter

\bmhead{Funding}
This work was supported by the National Natural Science Foundation of China (Grant No.~12631012).

\bmhead{Data availability}
No datasets were generated or analysed for this study.

\bmhead{Competing interests}
The authors declare no competing interests.

\bmhead{Disclosure}
During the preparation of this manuscript, the authors used generative AI-assisted tools for
language refinement, manuscript formatting, and reference organization. The authors reviewed
and revised all AI-assisted content and take full responsibility for the accuracy, analyses,
and conclusions presented in this work.

\end{document}